\documentclass[10pt,landscape,a4paper]{article}

\usepackage{multicol}
\usepackage{calc}
\usepackage{ifthen}
\usepackage[landscape]{geometry}
\usepackage{amsmath,amsthm,amsfonts,amssymb}
\usepackage{color,graphicx,overpic}
\usepackage{hyperref}
\usepackage[T1]{fontenc}
\usepackage{lmodern}         
\usepackage[utf8]{inputenc}  
\usepackage{cite}

\ifthenelse{\lengthtest { \paperwidth = 11in}}
    { \geometry{top=.5in,left=.5in,right=.5in,bottom=.5in} }
    {\ifthenelse{ \lengthtest{ \paperwidth = 297mm}}
        {\geometry{top=2cm,left=1cm,right=1cm,bottom=1cm} }
        {\geometry{top=2cm,left=1cm,right=1cm,bottom=1cm} }
    }

\usepackage{fancyhdr}
\makeatletter
\renewcommand{\section}{\@startsection{section}{1}{0mm}%
                                {-1ex plus -.5ex minus -.2ex}%
                                {0.5ex plus .2ex}
                                {\normalfont\Large\bfseries}}
\renewcommand{\subsection}{\@startsection{subsection}{2}{0mm}%
                                {-1explus -.5ex minus -.2ex}%
                                {0.5ex plus .2ex}%
                                {\normalfont\normalsize\bfseries}}
\renewcommand{\subsubsection}{\@startsection{subsubsection}{3}{0mm}%
                                {-1ex plus -.5ex minus -.2ex}%
                                {1ex plus .2ex}%
                                {\normalfont\small\bfseries}}
\makeatother

\def\BibTeX{{\rm B\kern-.05em{\sc i\kern-.025em b}\kern-.08em
    T\kern-.1667em\lower.7ex\hbox{E}\kern-.125emX}}

\begin{document}

\title{$3$-dimensional Continued Fraction Algorithms Cheat Sheets}
\author{Sébastien Labbé\footnote{
Université de Liège,
B\^at. B37 Institut de Math\'ematiques,
Grande Traverse 12,
4000 Li\`ege,
Belgium,
\texttt{slabbe@ulg.ac.be}.}}
\date{}

\maketitle

\thispagestyle{empty} 

\begin{multicols}{2}
\begin{abstract}
Multidimensional Continued Fraction Algorithms are generalizations of the
Euclid algorithm and find iteratively the gcd of two or more numbers. They are
defined as linear applications on some subcone of $\mathbb{R}^d$. We consider
multidimensional continued fraction algorithms that acts symmetrically on the
positive cone $\mathbb{R}^d_+$ for $d=3$. We include well-known and old ones
(Poincar\'e, Brun, Selmer, Fully Subtractive) and new ones
(Arnoux-Rauzy-Poincar\'e, Reverse, Cassaigne). 

For each algorithm, one page (called cheat sheet) gathers a handful of
informations most of them generated with the open source software Sage
\cite{sage} with the optional Sage package \texttt{slabbe-0.2.spkg}
\cite{labbe_slabbe_2015}. The information includes the $n$-cylinders, density
function of an absolutely continuous invariant measure, domain of the natural
extension, lyapunov exponents as well as data regarding combinatorics on
words, symbolic dynamics and digital geometry, that is, associated
substitutions, generated $S$-adic systems, factor complexity, discrepancy,
dual substitutions and generation of digital planes.

The document ends with a table of comparison of Lyapunov exponents and gives the
code allowing to reproduce any of the results or figures appearing in these
cheat sheets.
\end{abstract}

\columnbreak
\tableofcontents
\end{multicols}

\raggedright

\newpage

\begin{multicols}{3}
\setlength{\premulticols}{1pt}
\setlength{\postmulticols}{1pt}
\setlength{\multicolsep}{1pt}
\setlength{\columnsep}{2pt}
\newcommand{\note}[1]{\hfill\textrm{\textcolor{gray}{#1}}}
\newcommand{\args}[1]{\textit{\textcolor{blue}{#1}}}
\newcommand{\stdout}[1]{\textcolor{Sepia}{#1}}
\footnotesize

\section{Brun algorithm}
\subsection{Definition}
On $\Lambda=\mathbb{R}^3_+$, the map
\[
F (x_1,x_2,x_3) = (x'_1,x'_2,x'_3)
\]
is defined by
\[
    (x'_{\pi 1}, x'_{\pi 2}, x'_{\pi 3}) =
    (x_{\pi 1}, x_{\pi 2}, x_{\pi 3}-x_{\pi 2})
\]
where $\pi\in\mathcal{S}_3$ is the permutation of $\{1,2,3\}$ such that
$x_{\pi 1}<x_{\pi 2}<x_{\pi 3}$ \cite{MR0111735}.
\subsection{Matrix Definition}
The partition of the cone is
$\Lambda=\cup_{\pi\in\mathcal{S}_3}\Lambda_\pi$ where
\[
    \Lambda_\pi = \{(x_1,x_2,x_3)\in\Lambda\mid 
	x_{\pi 1}< x_{\pi 2}< x_{\pi 3}\}.
\]
The matrices are given by the rule
\[
    M(\mathbf{x}) = M_\pi
    \qquad\text{ if and only if }\qquad
    \mathbf{x}\in\Lambda_\pi.
\]
The map $F$ on $\Lambda$ and
the projective map $f$ on
$\Delta=\{\mathbf{x}\in\Lambda\mid\Vert\mathbf{x}\Vert_1=1\}$ are:
\[
    F(\mathbf{x}) = M(\mathbf{x})^{-1}\mathbf{x}
    \qquad\text{and}\qquad
    f(\mathbf{x}) = \frac{F(\mathbf{x})}{\Vert F(\mathbf{x})\Vert_1}.
\]

\subsection{Matrices}
\[
\begin{array}{lll}
M_{123}={\arraycolsep=2pt\left(\begin{array}{rrr}
1 & 0 & 0 \\
0 & 1 & 0 \\
0 & 1 & 1
\end{array}\right)}
&
M_{132}={\arraycolsep=2pt\left(\begin{array}{rrr}
1 & 0 & 0 \\
0 & 1 & 1 \\
0 & 0 & 1
\end{array}\right)}
&
M_{213}={\arraycolsep=2pt\left(\begin{array}{rrr}
1 & 0 & 0 \\
0 & 1 & 0 \\
1 & 0 & 1
\end{array}\right)}
\\
M_{231}={\arraycolsep=2pt\left(\begin{array}{rrr}
1 & 0 & 1 \\
0 & 1 & 0 \\
0 & 0 & 1
\end{array}\right)}
&
M_{312}={\arraycolsep=2pt\left(\begin{array}{rrr}
1 & 0 & 0 \\
1 & 1 & 0 \\
0 & 0 & 1
\end{array}\right)}
&
M_{321}={\arraycolsep=2pt\left(\begin{array}{rrr}
1 & 1 & 0 \\
0 & 1 & 0 \\
0 & 0 & 1
\end{array}\right)}
\\
\end{array}
\]
\subsection{Cylinders}
\includegraphics[width=0.300000000000000\linewidth]{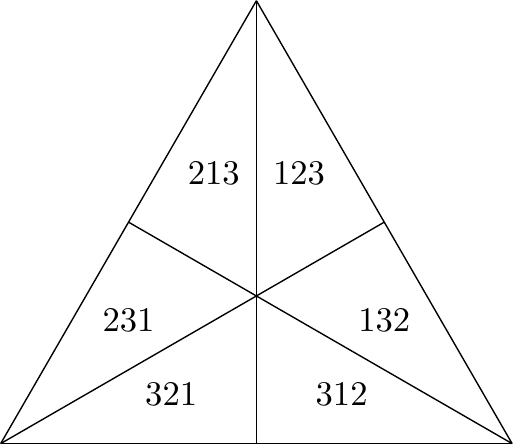}
\includegraphics[width=0.300000000000000\linewidth]{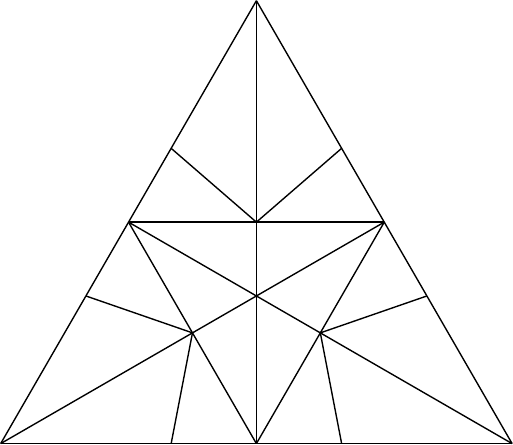}
\includegraphics[width=0.300000000000000\linewidth]{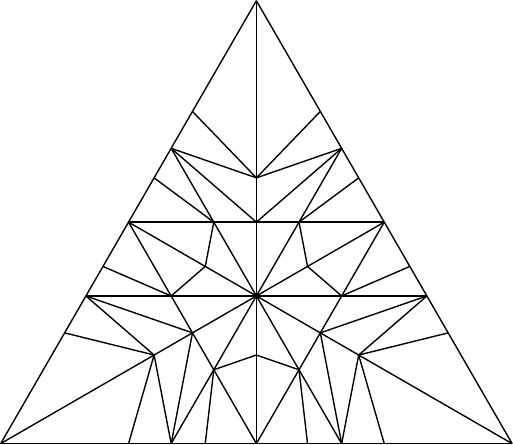}
\includegraphics[width=0.300000000000000\linewidth]{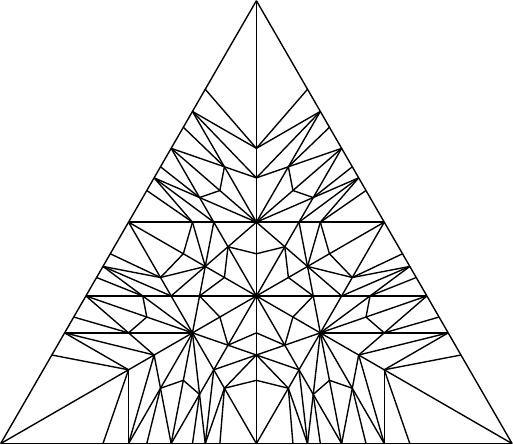}
\includegraphics[width=0.300000000000000\linewidth]{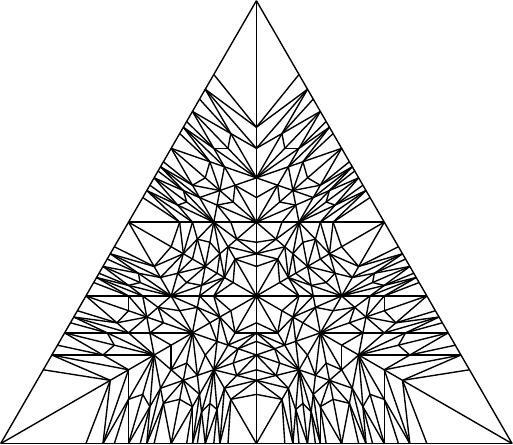}
\includegraphics[width=0.300000000000000\linewidth]{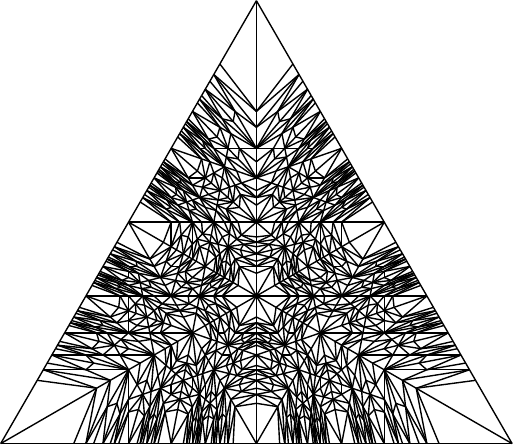}
\subsection{Density function}
The density function of the invariant measure of $f:\Delta\to\Delta$ for
the Brun algorithm is \cite{arnoux_symmetric_2015}:
\[
\frac{1}{2\,x_{\pi 2}(1-x_{\pi 2})(1-x_{\pi 1}-x_{\pi 2})}
\]
on the part $\mathbf{x}=(x_1,x_2,x_3)\in\Lambda_\pi\cap\Delta$.

\subsection{Invariant measure}
\includegraphics[width=0.800000000000000\linewidth]{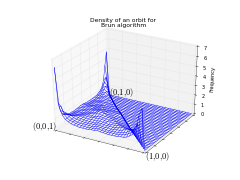}
\subsection{Natural extension}
Two sequences 
$(\mathbf{x}_{n+1})_{n\geq0}$ and 
$(\mathbf{a}_{n+1})_{n\geq0}$ defined such that
\[
    \mathbf{x}_{n+1}=M(\mathbf{x}_n)^{-1}\mathbf{x}_n
    \qquad\text{ and }\qquad
    \mathbf{a}_{n+1}=M(\mathbf{x}_n)^\top\mathbf{a}_n.
\]

\includegraphics[width=1\linewidth]{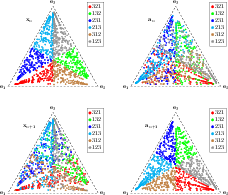}
\subsection{Lyapunov exponents}
(using 30 orbits of 
100000000 iterations each)\\
\begin{tabular}{lllll}
30 succesfull orbits & min & mean & max & std \\ \hline
$\theta_1$ & $0.3041$ & $0.3044$ & $0.3047$ & $0.00016$ \\
$\theta_2$ & $-0.11222$ & $-0.11213$ & $-0.11198$ & $0.000067$ \\
$1-\theta_2/\theta_1$ & $1.36821$ & $1.36834$ & $1.36845$ & $0.000063$ \\
\end{tabular}
\subsection{Substitutions}
\[
\begin{array}{lll}
\sigma_{123}=\left\{\begin{array}{l}
1 \mapsto 1\\
2 \mapsto 23\\
3 \mapsto 3
\end{array}\right.
&
\sigma_{132}=\left\{\begin{array}{l}
1 \mapsto 1\\
2 \mapsto 2\\
3 \mapsto 32
\end{array}\right.
&
\sigma_{213}=\left\{\begin{array}{l}
1 \mapsto 13\\
2 \mapsto 2\\
3 \mapsto 3
\end{array}\right.
\\
\sigma_{231}=\left\{\begin{array}{l}
1 \mapsto 1\\
2 \mapsto 2\\
3 \mapsto 31
\end{array}\right.
&
\sigma_{312}=\left\{\begin{array}{l}
1 \mapsto 12\\
2 \mapsto 2\\
3 \mapsto 3
\end{array}\right.
&
\sigma_{321}=\left\{\begin{array}{l}
1 \mapsto 1\\
2 \mapsto 21\\
3 \mapsto 3
\end{array}\right.
\\
\end{array}
\]
\subsection{$S$-adic word example}
Using vector $v=\left(1, e, \pi\right)$:
\begin{align*}
w &=
\sigma_{123}
\sigma_{312}
\sigma_{312}
\sigma_{321}
\sigma_{132}
\sigma_{123}
\sigma_{312}
\sigma_{231}
\sigma_{231}
\sigma_{213}
\cdots(1)\\
& = 1232323123233231232332312323123232312323...
\end{align*}
Factor Complexity of $w$ is 
$(p_w(n))_{0\leq n \leq 20} =$
\[
(1, 3, 5, 7, 9, 11, 13, 15, 17, 19, 22, 24, 26, 28, 30, 32, 34, 36, 38, 40, 42)
\]
\subsection{Discrepancy}
Discrepancy \cite{MR593979} for all 19701
 $S$-adic words with directions
$v\in\mathbb{N}^3_{>0}$
such that $v_1+v_2+v_3=200$:
\begin{center}
\includegraphics[width=0.600000000000000\linewidth]{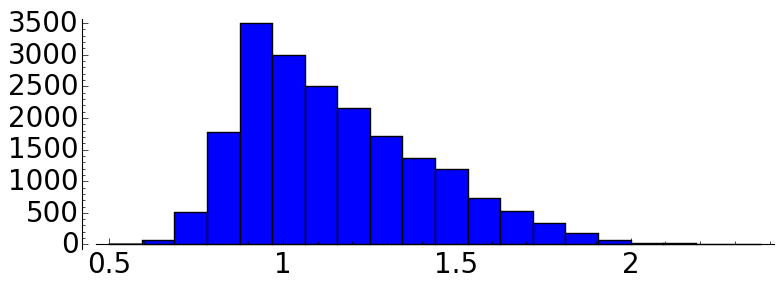}
\end{center}
\subsection{Dual substitutions}
\[
\begin{array}{lll}
\sigma^*_{123}=\left\{\begin{array}{l}
1 \mapsto 1\\
2 \mapsto 2\\
3 \mapsto 32
\end{array}\right.
&
\sigma^*_{132}=\left\{\begin{array}{l}
1 \mapsto 1\\
2 \mapsto 23\\
3 \mapsto 3
\end{array}\right.
&
\sigma^*_{213}=\left\{\begin{array}{l}
1 \mapsto 1\\
2 \mapsto 2\\
3 \mapsto 31
\end{array}\right.
\\
\sigma^*_{231}=\left\{\begin{array}{l}
1 \mapsto 13\\
2 \mapsto 2\\
3 \mapsto 3
\end{array}\right.
&
\sigma^*_{312}=\left\{\begin{array}{l}
1 \mapsto 1\\
2 \mapsto 21\\
3 \mapsto 3
\end{array}\right.
&
\sigma^*_{321}=\left\{\begin{array}{l}
1 \mapsto 12\\
2 \mapsto 2\\
3 \mapsto 3
\end{array}\right.
\\
\end{array}
\]
\subsection{E one star}
Using vector $v=\left(1, e, \pi\right)$, the 9-th 
iteration on the unit cube is:
\[
E_1^*(\sigma^*_{123})
E_1^*(\sigma^*_{312})
E_1^*(\sigma^*_{312})
E_1^*(\sigma^*_{321})
E_1^*(\sigma^*_{132})
\cdots
(\includegraphics[width=1em]{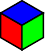})=
\]
\begin{center}
\includegraphics[height=3cm]{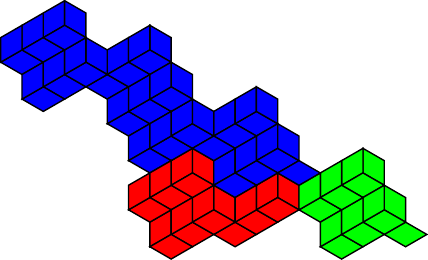}
\end{center}
\newpage
\section{Selmer algorithm}
\subsection{Definition}
On $\Lambda=\mathbb{R}^3_+$, the map
\[
F (x_1,x_2,x_3) = (x'_1,x'_2,x'_3)
\]
is defined by
\[
    (x'_{\pi 1}, x'_{\pi 2}, x'_{\pi 3}) =
    (x_{\pi 1}, x_{\pi 2}, x_{\pi 3}-x_{\pi 1})
\]
where $\pi\in\mathcal{S}_3$ is the permutation of $\{1,2,3\}$ such that
$x_{\pi 1}<x_{\pi 2}<x_{\pi 3}$
\cite{MR0130852}.
\subsection{Matrix Definition}
The partition of the cone is
$\Lambda=\cup_{\pi\in\mathcal{S}_3}\Lambda_\pi$ where
\[
    \Lambda_\pi = \{(x_1,x_2,x_3)\in\Lambda\mid 
	x_{\pi 1}< x_{\pi 2}< x_{\pi 3}\}.
\]
The matrices are given by the rule
\[
    M(\mathbf{x}) = M_\pi
    \qquad\text{ if and only if }\qquad
    \mathbf{x}\in\Lambda_\pi.
\]
The map $F$ on $\Lambda$ and
the projective map $f$ on
$\Delta=\{\mathbf{x}\in\Lambda\mid\Vert\mathbf{x}\Vert_1=1\}$ are:
\[
    F(\mathbf{x}) = M(\mathbf{x})^{-1}\mathbf{x}
    \qquad\text{and}\qquad
    f(\mathbf{x}) = \frac{F(\mathbf{x})}{\Vert F(\mathbf{x})\Vert_1}.
\]

\subsection{Matrices}
\[
\begin{array}{lll}
M_{123}={\arraycolsep=2pt\left(\begin{array}{rrr}
1 & 0 & 0 \\
0 & 1 & 0 \\
1 & 0 & 1
\end{array}\right)}
&
M_{132}={\arraycolsep=2pt\left(\begin{array}{rrr}
1 & 0 & 0 \\
1 & 1 & 0 \\
0 & 0 & 1
\end{array}\right)}
&
M_{213}={\arraycolsep=2pt\left(\begin{array}{rrr}
1 & 0 & 0 \\
0 & 1 & 0 \\
0 & 1 & 1
\end{array}\right)}
\\
M_{231}={\arraycolsep=2pt\left(\begin{array}{rrr}
1 & 1 & 0 \\
0 & 1 & 0 \\
0 & 0 & 1
\end{array}\right)}
&
M_{312}={\arraycolsep=2pt\left(\begin{array}{rrr}
1 & 0 & 0 \\
0 & 1 & 1 \\
0 & 0 & 1
\end{array}\right)}
&
M_{321}={\arraycolsep=2pt\left(\begin{array}{rrr}
1 & 0 & 1 \\
0 & 1 & 0 \\
0 & 0 & 1
\end{array}\right)}
\\
\end{array}
\]
\subsection{Cylinders}
\includegraphics[width=0.300000000000000\linewidth]{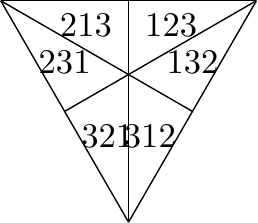}
\includegraphics[width=0.300000000000000\linewidth]{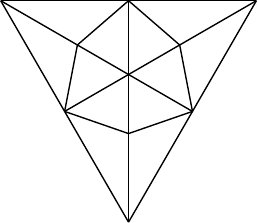}
\includegraphics[width=0.300000000000000\linewidth]{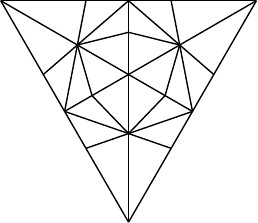}
\includegraphics[width=0.300000000000000\linewidth]{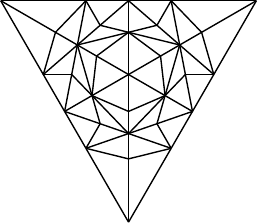}
\includegraphics[width=0.300000000000000\linewidth]{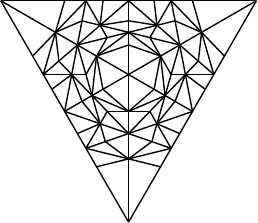}
\includegraphics[width=0.300000000000000\linewidth]{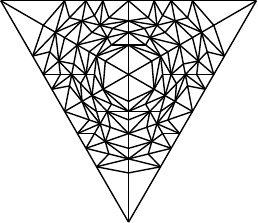}
\subsection{Density function}
The sorted version of $f$ admits a $\sigma$-finite invariant measure which is
absolutely continuous with respect to Lebesgue measure on the central part and
its density is known \cite{schweiger}.

\subsection{Invariant measure}
\includegraphics[width=0.800000000000000\linewidth]{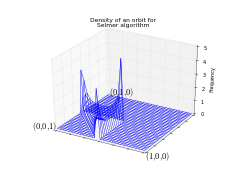}
\subsection{Natural extension}
Two sequences 
$(\mathbf{x}_{n+1})_{n\geq0}$ and 
$(\mathbf{a}_{n+1})_{n\geq0}$ defined such that
\[
    \mathbf{x}_{n+1}=M(\mathbf{x}_n)^{-1}\mathbf{x}_n
    \qquad\text{ and }\qquad
    \mathbf{a}_{n+1}=M(\mathbf{x}_n)^\top\mathbf{a}_n.
\]

\includegraphics[width=1\linewidth]{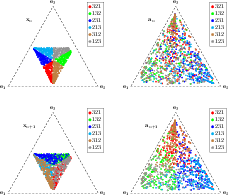}
\subsection{Lyapunov exponents}
(using 30 orbits of 
100000000 iterations each)\\
\begin{tabular}{lllll}
30 succesfull orbits & min & mean & max & std \\ \hline
$\theta_1$ & $0.1814$ & $0.1827$ & $0.1829$ & $0.00028$ \\
$\theta_2$ & $-0.0708$ & $-0.0707$ & $-0.0702$ & $0.00011$ \\
$1-\theta_2/\theta_1$ & $1.38697$ & $1.38711$ & $1.38726$ & $0.000084$ \\
\end{tabular}
\subsection{Substitutions}
\[
\begin{array}{lll}
\sigma_{123}=\left\{\begin{array}{l}
1 \mapsto 13\\
2 \mapsto 2\\
3 \mapsto 3
\end{array}\right.
&
\sigma_{132}=\left\{\begin{array}{l}
1 \mapsto 12\\
2 \mapsto 2\\
3 \mapsto 3
\end{array}\right.
&
\sigma_{213}=\left\{\begin{array}{l}
1 \mapsto 1\\
2 \mapsto 23\\
3 \mapsto 3
\end{array}\right.
\\
\sigma_{231}=\left\{\begin{array}{l}
1 \mapsto 1\\
2 \mapsto 21\\
3 \mapsto 3
\end{array}\right.
&
\sigma_{312}=\left\{\begin{array}{l}
1 \mapsto 1\\
2 \mapsto 2\\
3 \mapsto 32
\end{array}\right.
&
\sigma_{321}=\left\{\begin{array}{l}
1 \mapsto 1\\
2 \mapsto 2\\
3 \mapsto 31
\end{array}\right.
\\
\end{array}
\]
\subsection{$S$-adic word example}
Using vector $v=\left(1, e, \pi\right)$:
\begin{align*}
w &=
\sigma_{123}
\sigma_{132}
\sigma_{123}
\sigma_{132}
\sigma_{213}
\sigma_{321}
\sigma_{312}
\sigma_{231}
\sigma_{123}
\sigma_{312}
\cdots(1)\\
& = 1323231323223231323231323223231323213232...
\end{align*}
Factor Complexity of $w$ is 
$(p_w(n))_{0\leq n \leq 20} =$
\[
(1, 3, 7, 11, 16, 20, 24, 28, 32, 36, 40, 44, 48, 52, 56, 60, 64, 68, 72, 76, 80)
\]
\subsection{Discrepancy}
ValueError: On input=[198, 1, 1], algorithm Selmer loops on (1.0, 1.0, 0.0)
\subsection{Dual substitutions}
\[
\begin{array}{lll}
\sigma^*_{123}=\left\{\begin{array}{l}
1 \mapsto 1\\
2 \mapsto 2\\
3 \mapsto 31
\end{array}\right.
&
\sigma^*_{132}=\left\{\begin{array}{l}
1 \mapsto 1\\
2 \mapsto 21\\
3 \mapsto 3
\end{array}\right.
&
\sigma^*_{213}=\left\{\begin{array}{l}
1 \mapsto 1\\
2 \mapsto 2\\
3 \mapsto 32
\end{array}\right.
\\
\sigma^*_{231}=\left\{\begin{array}{l}
1 \mapsto 12\\
2 \mapsto 2\\
3 \mapsto 3
\end{array}\right.
&
\sigma^*_{312}=\left\{\begin{array}{l}
1 \mapsto 1\\
2 \mapsto 23\\
3 \mapsto 3
\end{array}\right.
&
\sigma^*_{321}=\left\{\begin{array}{l}
1 \mapsto 13\\
2 \mapsto 2\\
3 \mapsto 3
\end{array}\right.
\\
\end{array}
\]
\subsection{E one star}
Using vector $v=\left(1, e, \pi\right)$, the 13-th 
iteration on the unit cube is:
\[
E_1^*(\sigma^*_{123})
E_1^*(\sigma^*_{132})
E_1^*(\sigma^*_{123})
E_1^*(\sigma^*_{132})
E_1^*(\sigma^*_{213})
\cdots
(\includegraphics[width=1em]{cube.pdf})=
\]
\begin{center}
\includegraphics[height=3cm]{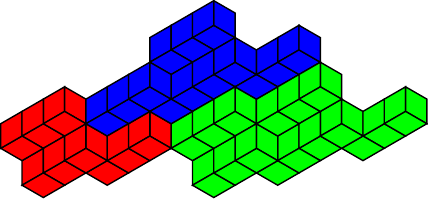}
\end{center}
\newpage
\section{Poincar\'e algorithm}
\subsection{Definition}
On $\Lambda=\mathbb{R}^3_+$, the map
\[
F (x_1,x_2,x_3) = (x'_1,x'_2,x'_3)
\]
is defined by
\[
    (x'_{\pi 1}, x'_{\pi 2}, x'_{\pi 3}) =
    (x_{\pi 1}, x_{\pi 2}-x_{\pi 1}, x_{\pi 3}-x_{\pi 2})
\]
where $\pi\in\mathcal{S}_3$ is the permutation of $\{1,2,3\}$ such that
$x_{\pi 1}<x_{\pi 2}<x_{\pi 3}$ \cite{MR1336331}.
\subsection{Matrix Definition}
The partition of the cone is
$\Lambda=\cup_{\pi\in\mathcal{S}_3}\Lambda_\pi$ where
\[
    \Lambda_\pi = \{(x_1,x_2,x_3)\in\Lambda\mid 
	x_{\pi 1}< x_{\pi 2}< x_{\pi 3}\}.
\]
The matrices are given by the rule
\[
    M(\mathbf{x}) = M_\pi
    \qquad\text{ if and only if }\qquad
    \mathbf{x}\in\Lambda_\pi.
\]
The map $F$ on $\Lambda$ and
the projective map $f$ on
$\Delta=\{\mathbf{x}\in\Lambda\mid\Vert\mathbf{x}\Vert_1=1\}$ are:
\[
    F(\mathbf{x}) = M(\mathbf{x})^{-1}\mathbf{x}
    \qquad\text{and}\qquad
    f(\mathbf{x}) = \frac{F(\mathbf{x})}{\Vert F(\mathbf{x})\Vert_1}.
\]

\subsection{Matrices}
\[
\begin{array}{lll}
M_{123}={\arraycolsep=2pt\left(\begin{array}{rrr}
1 & 0 & 0 \\
1 & 1 & 0 \\
1 & 1 & 1
\end{array}\right)}
&
M_{132}={\arraycolsep=2pt\left(\begin{array}{rrr}
1 & 0 & 0 \\
1 & 1 & 1 \\
1 & 0 & 1
\end{array}\right)}
&
M_{213}={\arraycolsep=2pt\left(\begin{array}{rrr}
1 & 1 & 0 \\
0 & 1 & 0 \\
1 & 1 & 1
\end{array}\right)}
\\
M_{231}={\arraycolsep=2pt\left(\begin{array}{rrr}
1 & 1 & 1 \\
0 & 1 & 0 \\
0 & 1 & 1
\end{array}\right)}
&
M_{312}={\arraycolsep=2pt\left(\begin{array}{rrr}
1 & 0 & 1 \\
1 & 1 & 1 \\
0 & 0 & 1
\end{array}\right)}
&
M_{321}={\arraycolsep=2pt\left(\begin{array}{rrr}
1 & 1 & 1 \\
0 & 1 & 1 \\
0 & 0 & 1
\end{array}\right)}
\\
\end{array}
\]
\subsection{Cylinders}
\includegraphics[width=0.300000000000000\linewidth]{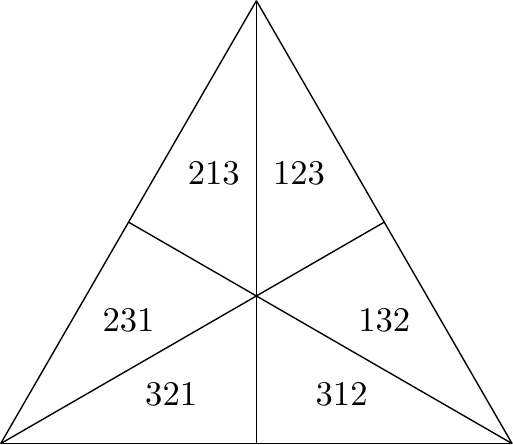}
\includegraphics[width=0.300000000000000\linewidth]{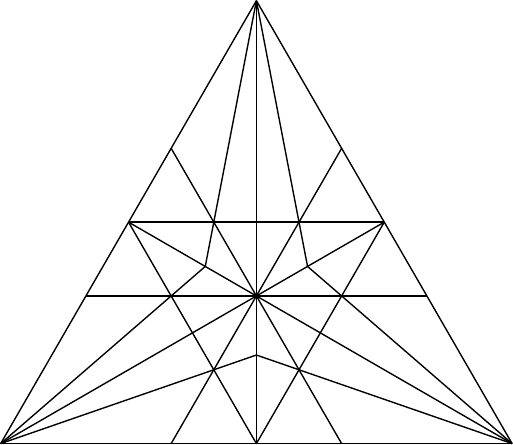}
\includegraphics[width=0.300000000000000\linewidth]{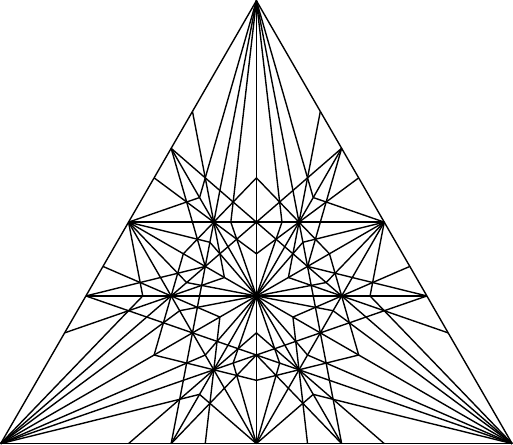}
\includegraphics[width=0.300000000000000\linewidth]{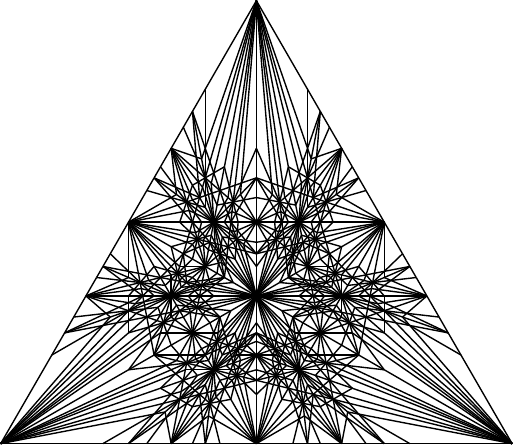}
\subsection{Density function}
The sorted version of $f$ admits a $\sigma$-finite invariant measure which is
absolutely continuous with respect to Lebesgue measure and its density is
known \cite{schweiger,MR1336331}.


\subsection{Invariant measure}
\includegraphics[width=0.800000000000000\linewidth]{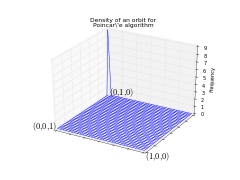}
\subsection{Natural extension}
Two sequences 
$(\mathbf{x}_{n+1})_{n\geq0}$ and 
$(\mathbf{a}_{n+1})_{n\geq0}$ defined such that
\[
    \mathbf{x}_{n+1}=M(\mathbf{x}_n)^{-1}\mathbf{x}_n
    \qquad\text{ and }\qquad
    \mathbf{a}_{n+1}=M(\mathbf{x}_n)^\top\mathbf{a}_n.
\]

\includegraphics[width=1\linewidth]{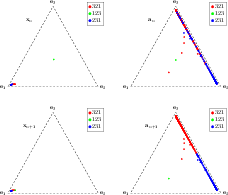}
\subsection{Lyapunov exponents}
(using 30 orbits of 
100000000 iterations each)\\
\begin{tabular}{lllll}
20 succesfull orbits & min & mean & max & std \\ \hline
$\theta_1$ & $6.1 \times 10^{-9}$ & $1.0 \times 10^{-7}$ & $2.9 \times 10^{-7}$ & $8.8 \times 10^{-8}$ \\
$\theta_2$ & $-8.2 \times 10^{-8}$ & $6.6 \times 10^{-7}$ & $9.5 \times 10^{-7}$ & $4.0 \times 10^{-7}$ \\
$1-\theta_2/\theta_1$ & $-150.$ & $-25.$ & $1.4$ & $40.$ \\
\end{tabular}
\subsection{Substitutions}
\[
\begin{array}{lll}
\sigma_{123}=\left\{\begin{array}{l}
1 \mapsto 123\\
2 \mapsto 23\\
3 \mapsto 3
\end{array}\right.
&
\sigma_{132}=\left\{\begin{array}{l}
1 \mapsto 132\\
2 \mapsto 2\\
3 \mapsto 32
\end{array}\right.
&
\sigma_{213}=\left\{\begin{array}{l}
1 \mapsto 13\\
2 \mapsto 213\\
3 \mapsto 3
\end{array}\right.
\\
\sigma_{231}=\left\{\begin{array}{l}
1 \mapsto 1\\
2 \mapsto 231\\
3 \mapsto 31
\end{array}\right.
&
\sigma_{312}=\left\{\begin{array}{l}
1 \mapsto 12\\
2 \mapsto 2\\
3 \mapsto 312
\end{array}\right.
&
\sigma_{321}=\left\{\begin{array}{l}
1 \mapsto 1\\
2 \mapsto 21\\
3 \mapsto 321
\end{array}\right.
\\
\end{array}
\]
\subsection{$S$-adic word example}
Using vector $v=\left(1, e, \pi\right)$:
\begin{align*}
w &=
\sigma_{123}
\sigma_{312}
\sigma_{312}
\sigma_{213}
\sigma_{123}
\sigma_{132}
\sigma_{213}
\sigma_{213}
\sigma_{213}
\sigma_{213}
\cdots(1)\\
& = 1232323312323123232323123232331232312323...
\end{align*}
Factor Complexity of $w$ is 
$(p_w(n))_{0\leq n \leq 20} =$
\[
(1, 3, 5, 7, 9, 11, 14, 17, 19, 21, 23, 25, 27, 29, 31, 33, 35, 37, 39, 41, 43)
\]
\subsection{Discrepancy}
Discrepancy \cite{MR593979} for all 19701
 $S$-adic words with directions
$v\in\mathbb{N}^3_{>0}$
such that $v_1+v_2+v_3=200$:
\begin{center}
\includegraphics[width=0.600000000000000\linewidth]{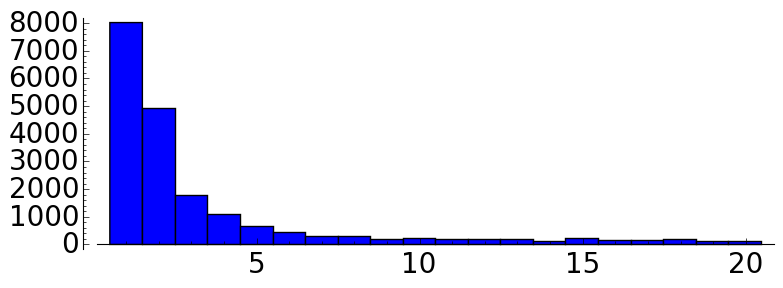}
\end{center}
\subsection{Dual substitutions}
\[
\begin{array}{lll}
\sigma^*_{123}=\left\{\begin{array}{l}
1 \mapsto 1\\
2 \mapsto 21\\
3 \mapsto 321
\end{array}\right.
&
\sigma^*_{132}=\left\{\begin{array}{l}
1 \mapsto 1\\
2 \mapsto 231\\
3 \mapsto 31
\end{array}\right.
&
\sigma^*_{213}=\left\{\begin{array}{l}
1 \mapsto 12\\
2 \mapsto 2\\
3 \mapsto 312
\end{array}\right.
\\
\sigma^*_{231}=\left\{\begin{array}{l}
1 \mapsto 132\\
2 \mapsto 2\\
3 \mapsto 32
\end{array}\right.
&
\sigma^*_{312}=\left\{\begin{array}{l}
1 \mapsto 13\\
2 \mapsto 213\\
3 \mapsto 3
\end{array}\right.
&
\sigma^*_{321}=\left\{\begin{array}{l}
1 \mapsto 123\\
2 \mapsto 23\\
3 \mapsto 3
\end{array}\right.
\\
\end{array}
\]
\subsection{E one star}
Using vector $v=\left(1, e, \pi\right)$, the 5-th 
iteration on the unit cube is:
\[
E_1^*(\sigma^*_{123})
E_1^*(\sigma^*_{312})
E_1^*(\sigma^*_{312})
E_1^*(\sigma^*_{213})
E_1^*(\sigma^*_{123})
(\includegraphics[width=1em]{cube.pdf})=
\]
\begin{center}
\includegraphics[height=3cm]{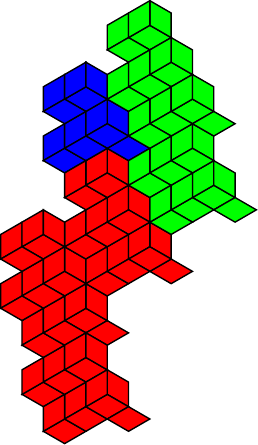}
\end{center}
\newpage
\section{Fully Subtractive algorithm}
\subsection{Definition}
On $\Lambda=\mathbb{R}^3_+$, the map
\[
F (x_1,x_2,x_3) = (x'_1,x'_2,x'_3)
\]
is defined by
\[
    (x'_{\pi 1}, x'_{\pi 2}, x'_{\pi 3}) =
    (x_{\pi 1}, x_{\pi 2}-x_{\pi 1}, x_{\pi 3}-x_{\pi 1})
\]
where $\pi\in\mathcal{S}_3$ is the permutation of $\{1,2,3\}$ such that
$x_{\pi 1}<x_{\pi 2}<x_{\pi 3}$
\cite{schweiger}.
\subsection{Matrix Definition}
The partition of the cone is
$\Lambda=\cup_{i\in\{1,2,3\}}\Lambda_i$ where
\[
    \Lambda_i = \left\{(x_1,x_2,x_3)\in\Lambda\mid 
    x_i = \min\{x_1,x_2,x_3\}\right\}.
\]
The matrices are given by the rule
\[
    M(\mathbf{x}) = M_i
    \qquad\text{ if and only if }\qquad
    \mathbf{x}\in\Lambda_i.
\]
The map $F$ on $\Lambda$ and
the projective map $f$ on
$\Delta=\{\mathbf{x}\in\Lambda\mid\Vert\mathbf{x}\Vert_1=1\}$ are:
\[
    F(\mathbf{x}) = M(\mathbf{x})^{-1}\mathbf{x}
    \qquad\text{and}\qquad
    f(\mathbf{x}) = \frac{F(\mathbf{x})}{\Vert F(\mathbf{x})\Vert_1}.
\]

\subsection{Matrices}
\[
\begin{array}{lll}
M_{1}={\arraycolsep=2pt\left(\begin{array}{rrr}
1 & 0 & 0 \\
1 & 1 & 0 \\
1 & 0 & 1
\end{array}\right)}
&
M_{2}={\arraycolsep=2pt\left(\begin{array}{rrr}
1 & 1 & 0 \\
0 & 1 & 0 \\
0 & 1 & 1
\end{array}\right)}
&
M_{3}={\arraycolsep=2pt\left(\begin{array}{rrr}
1 & 0 & 1 \\
0 & 1 & 1 \\
0 & 0 & 1
\end{array}\right)}
\\
\end{array}
\]
\subsection{Cylinders}
\includegraphics[width=0.300000000000000\linewidth]{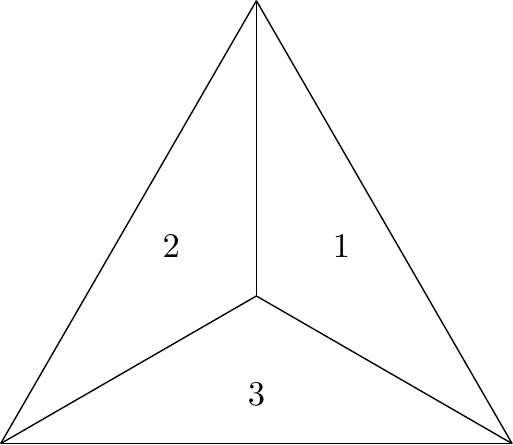}
\includegraphics[width=0.300000000000000\linewidth]{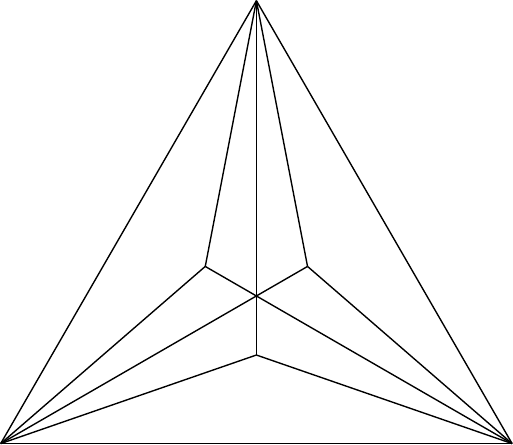}
\includegraphics[width=0.300000000000000\linewidth]{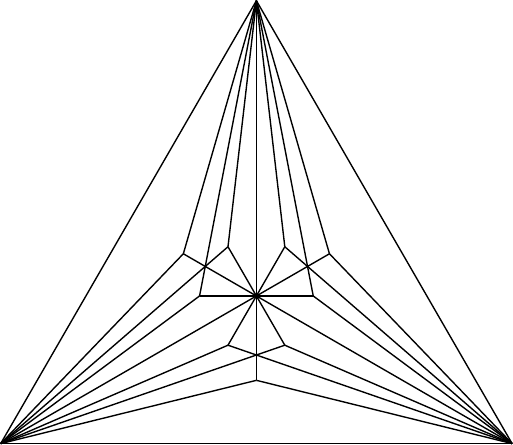}
\includegraphics[width=0.300000000000000\linewidth]{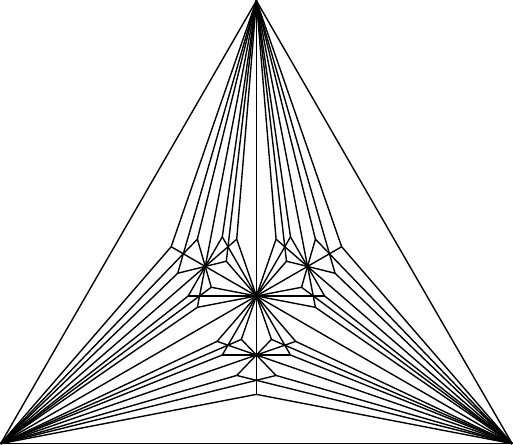}
\includegraphics[width=0.300000000000000\linewidth]{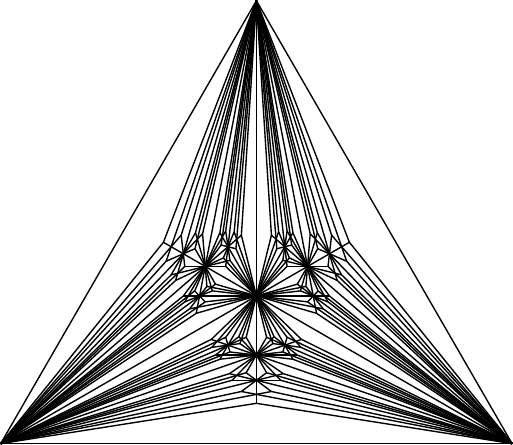}
\includegraphics[width=0.300000000000000\linewidth]{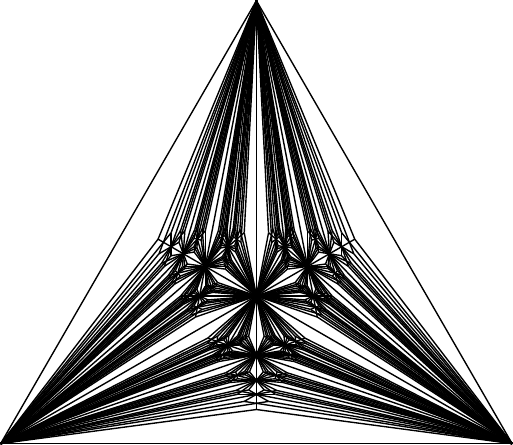}
\subsection{Density function}
The sorted version of $f$ admits a $\sigma$-finite invariant measure which is
absolutely continuous with respect to Lebesgue measure and its density is
known \cite{schweiger}.

\subsection{Invariant measure}
\includegraphics[width=0.800000000000000\linewidth]{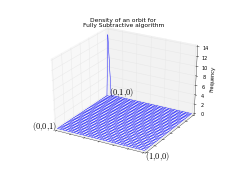}
\subsection{Natural extension}
Two sequences 
$(\mathbf{x}_{n+1})_{n\geq0}$ and 
$(\mathbf{a}_{n+1})_{n\geq0}$ defined such that
\[
    \mathbf{x}_{n+1}=M(\mathbf{x}_n)^{-1}\mathbf{x}_n
    \qquad\text{ and }\qquad
    \mathbf{a}_{n+1}=M(\mathbf{x}_n)^\top\mathbf{a}_n.
\]

\includegraphics[width=1\linewidth]{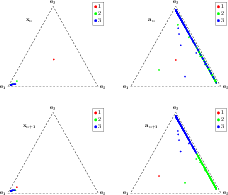}
\subsection{Lyapunov exponents}
(using 30 orbits of 
100000000 iterations each)\\
\begin{tabular}{lllll}
24 succesfull orbits & min & mean & max & std \\ \hline
$\theta_1$ & $2.6 \times 10^{-9}$ & $2.8 \times 10^{-8}$ & $6.1 \times 10^{-8}$ & $1.7 \times 10^{-8}$ \\
$\theta_2$ & $8.5 \times 10^{-7}$ & $9.1 \times 10^{-7}$ & $9.8 \times 10^{-7}$ & $2.4 \times 10^{-8}$ \\
$1-\theta_2/\theta_1$ & $-360.$ & $-63.$ & $-14.$ & $79.$ \\
\end{tabular}
\subsection{Substitutions}
\[
\begin{array}{lll}
\sigma_{1}=\left\{\begin{array}{l}
1 \mapsto 123\\
2 \mapsto 2\\
3 \mapsto 3
\end{array}\right.
&
\sigma_{2}=\left\{\begin{array}{l}
1 \mapsto 1\\
2 \mapsto 231\\
3 \mapsto 3
\end{array}\right.
&
\sigma_{3}=\left\{\begin{array}{l}
1 \mapsto 1\\
2 \mapsto 2\\
3 \mapsto 312
\end{array}\right.
\\
\end{array}
\]
\subsection{$S$-adic word example}
Using vector $v=\left(1, e, \pi\right)$:
\begin{align*}
w &=
\sigma_{1}
\sigma_{1}
\sigma_{2}
\sigma_{1}
\sigma_{3}
\sigma_{1}
\sigma_{3}
\sigma_{3}
\sigma_{3}
\sigma_{3}
\cdots(1)\\
& = 1232323123233231232331232323123233231232...
\end{align*}
Factor Complexity of $w$ is 
$(p_w(n))_{0\leq n \leq 20} =$
\[
(1, 3, 5, 8, 11, 14, 16, 18, 19, 20, 21, 21, 21, 21, 21, 21, 21, 21, 21, 21, 21)
\]
\subsection{Discrepancy}
ValueError: On input=[198, 1, 1], algorithm Fully Subtractive loops on (197.0, 1.0, 0.0)
\subsection{Dual substitutions}
\[
\begin{array}{lll}
\sigma^*_{1}=\left\{\begin{array}{l}
1 \mapsto 1\\
2 \mapsto 21\\
3 \mapsto 31
\end{array}\right.
&
\sigma^*_{2}=\left\{\begin{array}{l}
1 \mapsto 12\\
2 \mapsto 2\\
3 \mapsto 32
\end{array}\right.
&
\sigma^*_{3}=\left\{\begin{array}{l}
1 \mapsto 13\\
2 \mapsto 23\\
3 \mapsto 3
\end{array}\right.
\\
\end{array}
\]
\subsection{E one star}
Using vector $v=\left(1, e, \pi\right)$, the 7-th 
iteration on the unit cube is:
\[
E_1^*(\sigma^*_{1})
E_1^*(\sigma^*_{1})
E_1^*(\sigma^*_{2})
E_1^*(\sigma^*_{1})
E_1^*(\sigma^*_{3})
\cdots
(\includegraphics[width=1em]{cube.pdf})=
\]
\begin{center}
\includegraphics[height=3cm]{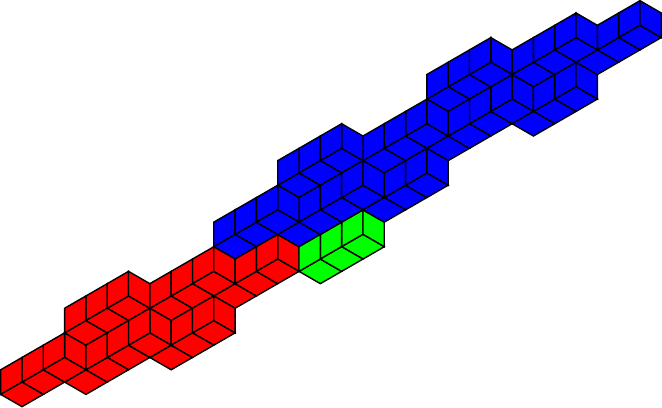}
\end{center}
\newpage
\section{Arnoux-Rauzy-Poincar\'e algorithm}
\subsection{Definition}
On $\Lambda=\mathbb{R}^3_+$, the map
\[
F (x_1,x_2,x_3) = (x'_1,x'_2,x'_3)
\]
is defined by
\[
(x'_{\pi 1}, x'_{\pi 2}, x'_{\pi 3}) =
\begin{cases}
    (x_{\pi 1}, x_{\pi 2}, x_{\pi 3}-x_{\pi 1}-x_{\pi 2}) &\mbox{if }
	x_{\pi 3}>x_{\pi 1}+x_{\pi 2}\\
    (x_{\pi 1}, x_{\pi 2}-x_{\pi 1}, x_{\pi 3}-x_{\pi 2}) &\mbox{otherwise.}
\end{cases}
\]
where $\pi\in\mathcal{S}_3$ is the permutation of $\{1,2,3\}$ such that
$x_{\pi 1}<x_{\pi 2}<x_{\pi 3}$ \cite{2015_berthe_factor}.
\subsection{Matrix Definition}
The subcones are
\begin{align*}
	\Lambda_i &= \{(x_1,x_2,x_3)\in\Lambda\mid 
	2x_i > x_1+x_2+x_3\},
	&i\in\{1,2,3\},\\
    \Lambda_\pi &= \{(x_1,x_2,x_3)\in\Lambda\mid 
	x_{\pi 1}< x_{\pi 2}< x_{\pi 3}\},
	&\pi\in\mathcal{S}_3.
\end{align*}
The matrices are given by the rule
\[
    M(\mathbf{x}) =
\begin{cases}
    M_i    & \text{ if } \mathbf{x} \in \Lambda_i,\\
    M_\pi  & \text{ else if } \mathbf{x} \in \Lambda_\pi.
\end{cases}
\]
The map $F$ on $\Lambda$ and
the projective map $f$ on
$\Delta=\{\mathbf{x}\in\Lambda\mid\Vert\mathbf{x}\Vert_1=1\}$ are:
\[
    F(\mathbf{x}) = M(\mathbf{x})^{-1}\mathbf{x}
    \qquad\text{and}\qquad
    f(\mathbf{x}) = \frac{F(\mathbf{x})}{\Vert F(\mathbf{x})\Vert_1}.
\]

\subsection{Matrices}
\[
\begin{array}{lll}
M_{1}={\arraycolsep=2pt\left(\begin{array}{rrr}
1 & 1 & 1 \\
0 & 1 & 0 \\
0 & 0 & 1
\end{array}\right)}
&
M_{2}={\arraycolsep=2pt\left(\begin{array}{rrr}
1 & 0 & 0 \\
1 & 1 & 1 \\
0 & 0 & 1
\end{array}\right)}
&
M_{3}={\arraycolsep=2pt\left(\begin{array}{rrr}
1 & 0 & 0 \\
0 & 1 & 0 \\
1 & 1 & 1
\end{array}\right)}
\\
M_{123}={\arraycolsep=2pt\left(\begin{array}{rrr}
1 & 0 & 0 \\
1 & 1 & 0 \\
1 & 1 & 1
\end{array}\right)}
&
M_{132}={\arraycolsep=2pt\left(\begin{array}{rrr}
1 & 0 & 0 \\
1 & 1 & 1 \\
1 & 0 & 1
\end{array}\right)}
&
M_{213}={\arraycolsep=2pt\left(\begin{array}{rrr}
1 & 1 & 0 \\
0 & 1 & 0 \\
1 & 1 & 1
\end{array}\right)}
\\
M_{231}={\arraycolsep=2pt\left(\begin{array}{rrr}
1 & 1 & 1 \\
0 & 1 & 0 \\
0 & 1 & 1
\end{array}\right)}
&
M_{312}={\arraycolsep=2pt\left(\begin{array}{rrr}
1 & 0 & 1 \\
1 & 1 & 1 \\
0 & 0 & 1
\end{array}\right)}
&
M_{321}={\arraycolsep=2pt\left(\begin{array}{rrr}
1 & 1 & 1 \\
0 & 1 & 1 \\
0 & 0 & 1
\end{array}\right)}
\\
\end{array}
\]
\subsection{Cylinders}
\includegraphics[width=0.300000000000000\linewidth]{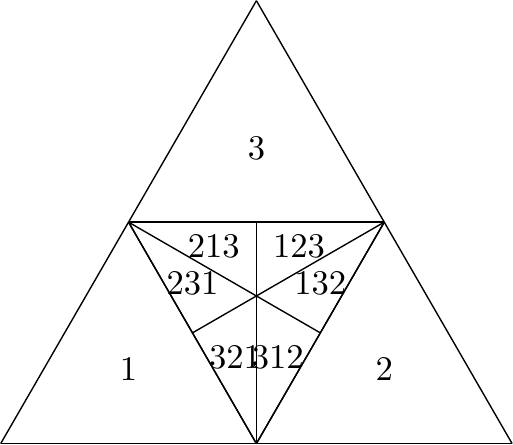}
\includegraphics[width=0.300000000000000\linewidth]{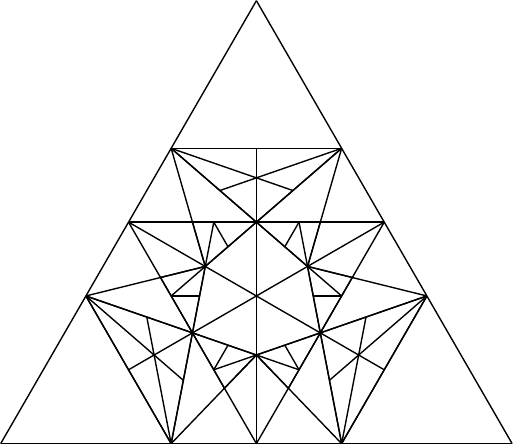}
\includegraphics[width=0.300000000000000\linewidth]{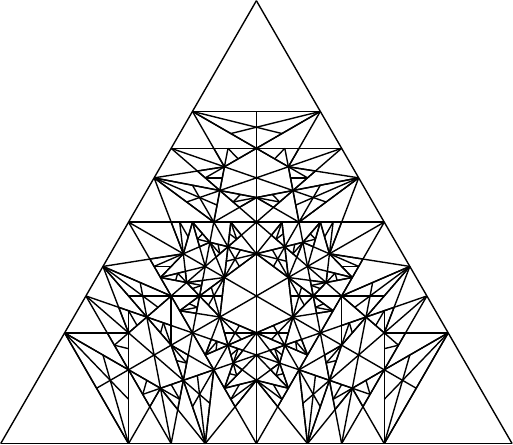}
\subsection{Density function}
The density of the absolutely continuous invariant measure is unknown
\cite{arnoux_symmetric_2015}.

\subsection{Invariant measure}
\includegraphics[width=0.800000000000000\linewidth]{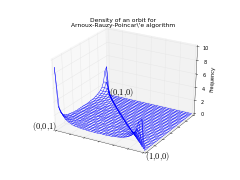}
\subsection{Natural extension}
Two sequences 
$(\mathbf{x}_{n+1})_{n\geq0}$ and 
$(\mathbf{a}_{n+1})_{n\geq0}$ defined such that
\[
    \mathbf{x}_{n+1}=M(\mathbf{x}_n)^{-1}\mathbf{x}_n
    \qquad\text{ and }\qquad
    \mathbf{a}_{n+1}=M(\mathbf{x}_n)^\top\mathbf{a}_n.
\]

\includegraphics[width=1\linewidth]{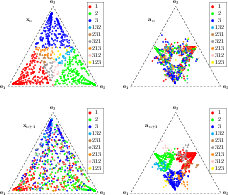}
\subsection{Lyapunov exponents}
(using 30 orbits of 
100000000 iterations each)\\
\begin{tabular}{lllll}
30 succesfull orbits & min & mean & max & std \\ \hline
$\theta_1$ & $0.4424$ & $0.4428$ & $0.4433$ & $0.00023$ \\
$\theta_2$ & $-0.17238$ & $-0.17218$ & $-0.17197$ & $0.000093$ \\
$1-\theta_2/\theta_1$ & $1.38873$ & $1.38880$ & $1.38890$ & $0.000049$ \\
\end{tabular}
\subsection{Substitutions}
\[
\begin{array}{lll}
\sigma_{1}=\left\{\begin{array}{l}
1 \mapsto 1\\
2 \mapsto 21\\
3 \mapsto 31
\end{array}\right.
&
\sigma_{2}=\left\{\begin{array}{l}
1 \mapsto 12\\
2 \mapsto 2\\
3 \mapsto 32
\end{array}\right.
&
\sigma_{3}=\left\{\begin{array}{l}
1 \mapsto 13\\
2 \mapsto 23\\
3 \mapsto 3
\end{array}\right.
\\
\sigma_{123}=\left\{\begin{array}{l}
1 \mapsto 123\\
2 \mapsto 23\\
3 \mapsto 3
\end{array}\right.
&
\sigma_{132}=\left\{\begin{array}{l}
1 \mapsto 132\\
2 \mapsto 2\\
3 \mapsto 32
\end{array}\right.
&
\sigma_{213}=\left\{\begin{array}{l}
1 \mapsto 13\\
2 \mapsto 213\\
3 \mapsto 3
\end{array}\right.
\\
\sigma_{231}=\left\{\begin{array}{l}
1 \mapsto 1\\
2 \mapsto 231\\
3 \mapsto 31
\end{array}\right.
&
\sigma_{312}=\left\{\begin{array}{l}
1 \mapsto 12\\
2 \mapsto 2\\
3 \mapsto 312
\end{array}\right.
&
\sigma_{321}=\left\{\begin{array}{l}
1 \mapsto 1\\
2 \mapsto 21\\
3 \mapsto 321
\end{array}\right.
\\
\end{array}
\]
\subsection{$S$-adic word example}
Using vector $v=\left(1, e, \pi\right)$:
\begin{align*}
w &=
\sigma_{123}
\sigma_{2}
\sigma_{1}
\sigma_{123}
\sigma_{1}
\sigma_{231}
\sigma_{3}
\sigma_{3}
\sigma_{3}
\sigma_{3}
\cdots(1)\\
& = 1232323123233231232332312323123232312323...
\end{align*}
Factor Complexity of $w$ is 
$(p_w(n))_{0\leq n \leq 20} =$
\[
(1, 3, 5, 7, 9, 11, 13, 15, 17, 19, 22, 24, 26, 28, 30, 32, 34, 36, 38, 40, 42)
\]
\subsection{Discrepancy}
Discrepancy \cite{MR593979} for all 19701
 $S$-adic words with directions
$v\in\mathbb{N}^3_{>0}$
such that $v_1+v_2+v_3=200$:
\begin{center}
\includegraphics[width=0.600000000000000\linewidth]{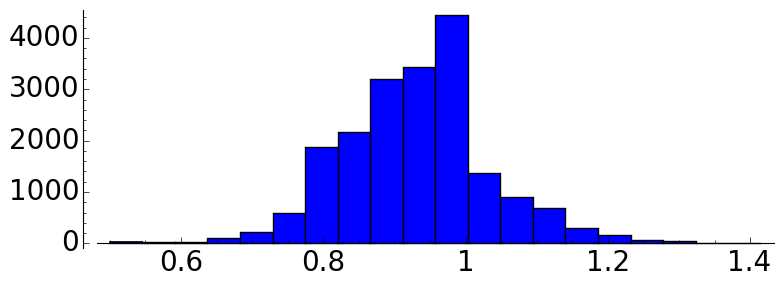}
\end{center}
\subsection{Dual substitutions}
\[
\begin{array}{lll}
\sigma^*_{1}=\left\{\begin{array}{l}
1 \mapsto 123\\
2 \mapsto 2\\
3 \mapsto 3
\end{array}\right.
&
\sigma^*_{2}=\left\{\begin{array}{l}
1 \mapsto 1\\
2 \mapsto 231\\
3 \mapsto 3
\end{array}\right.
&
\sigma^*_{3}=\left\{\begin{array}{l}
1 \mapsto 1\\
2 \mapsto 2\\
3 \mapsto 312
\end{array}\right.
\\
\sigma^*_{123}=\left\{\begin{array}{l}
1 \mapsto 1\\
2 \mapsto 21\\
3 \mapsto 321
\end{array}\right.
&
\sigma^*_{132}=\left\{\begin{array}{l}
1 \mapsto 1\\
2 \mapsto 231\\
3 \mapsto 31
\end{array}\right.
&
\sigma^*_{213}=\left\{\begin{array}{l}
1 \mapsto 12\\
2 \mapsto 2\\
3 \mapsto 312
\end{array}\right.
\\
\sigma^*_{231}=\left\{\begin{array}{l}
1 \mapsto 132\\
2 \mapsto 2\\
3 \mapsto 32
\end{array}\right.
&
\sigma^*_{312}=\left\{\begin{array}{l}
1 \mapsto 13\\
2 \mapsto 213\\
3 \mapsto 3
\end{array}\right.
&
\sigma^*_{321}=\left\{\begin{array}{l}
1 \mapsto 123\\
2 \mapsto 23\\
3 \mapsto 3
\end{array}\right.
\\
\end{array}
\]
\subsection{E one star}
Using vector $v=\left(1, e, \pi\right)$, the 5-th 
iteration on the unit cube is:
\[
E_1^*(\sigma^*_{123})
E_1^*(\sigma^*_{2})
E_1^*(\sigma^*_{1})
E_1^*(\sigma^*_{123})
E_1^*(\sigma^*_{1})
(\includegraphics[width=1em]{cube.pdf})=
\]
\begin{center}
\includegraphics[height=3cm]{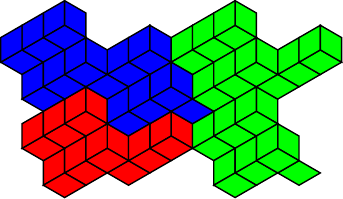}
\end{center}
\newpage
\section{Reverse algorithm}
\subsection{Definition}
On $\Lambda=\mathbb{R}^3_+$, the map
\[
F (x_1,x_2,x_3) = (x'_1,x'_2,x'_3)
\]
is defined by
\[
\left(\begin{array}{r}
	x'_{\pi 1} \\
	x'_{\pi 2} \\
	x'_{\pi 3}
\end{array}\right) =
\begin{cases}
\left(\begin{array}{l}
    x_{\pi 1}\\
    x_{\pi 2}\\
    x_{\pi 3}-x_{\pi 1}-x_{\pi 2}
\end{array}\right)
    &\mbox{if } x_{\pi 3}>x_{\pi 1}+x_{\pi 2}\\
\frac{1}{2}
\left(\begin{array}{r}
     -x_{\pi 1}+x_{\pi 2}+x_{\pi 3}\\
      x_{\pi 1}-x_{\pi 2}+x_{\pi 3}\\
      x_{\pi 1}+x_{\pi 2}-x_{\pi 3}
\end{array}\right)
     &\mbox{otherwise.}
\end{cases}
\]
where $\pi\in\mathcal{S}_3$ is the permutation of $\{1,2,3\}$ such that
$x_{\pi 1}<x_{\pi 2}<x_{\pi 3}$
\cite{arnoux_symmetric_2015}.
\subsection{Matrix Definition}
The subcones are
\begin{align*}
	\Lambda_i &= \{(x_1,x_2,x_3)\in\Lambda\mid 
	2x_i > x_1+x_2+x_3\},
	&i\in\{1,2,3\},\\
	\Lambda_4 &= \Lambda\setminus(\Lambda_1\cup\Lambda_2\cup\Lambda_3)
\end{align*}
The matrices are given by the rule
\[
    M(\mathbf{x}) = M_i
    \qquad\text{ if and only if }\qquad
    \mathbf{x}\in\Lambda_i.
\]
The map $F$ on $\Lambda$ and
the projective map $f$ on
$\Delta=\{\mathbf{x}\in\Lambda\mid\Vert\mathbf{x}\Vert_1=1\}$ are:
\[
    F(\mathbf{x}) = M(\mathbf{x})^{-1}\mathbf{x}
    \qquad\text{and}\qquad
    f(\mathbf{x}) = \frac{F(\mathbf{x})}{\Vert F(\mathbf{x})\Vert_1}.
\]

\subsection{Matrices}
\[
\begin{array}{lll}
M_{1}={\arraycolsep=2pt\left(\begin{array}{rrr}
1 & 1 & 1 \\
0 & 1 & 0 \\
0 & 0 & 1
\end{array}\right)}
&
M_{2}={\arraycolsep=2pt\left(\begin{array}{rrr}
1 & 0 & 0 \\
1 & 1 & 1 \\
0 & 0 & 1
\end{array}\right)}
&
M_{3}={\arraycolsep=2pt\left(\begin{array}{rrr}
1 & 0 & 0 \\
0 & 1 & 0 \\
1 & 1 & 1
\end{array}\right)}
\\
M_{4}={\arraycolsep=2pt\left(\begin{array}{rrr}
0 & 1 & 1 \\
1 & 0 & 1 \\
1 & 1 & 0
\end{array}\right)}
&
\end{array}
\]
\subsection{Cylinders}
\includegraphics[width=0.300000000000000\linewidth]{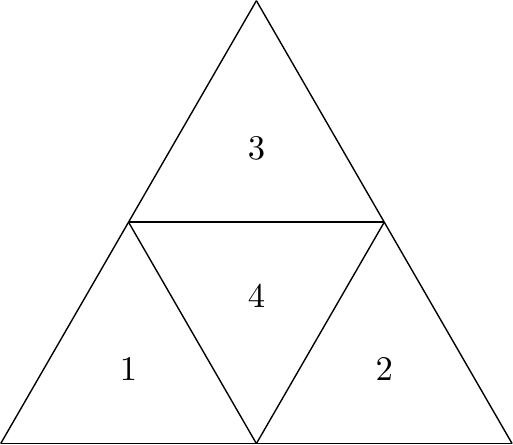}
\includegraphics[width=0.300000000000000\linewidth]{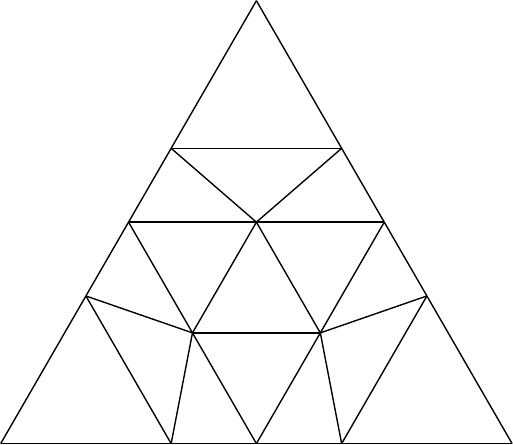}
\includegraphics[width=0.300000000000000\linewidth]{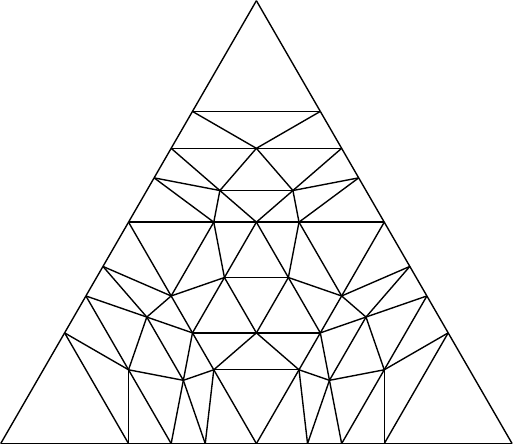}
\includegraphics[width=0.300000000000000\linewidth]{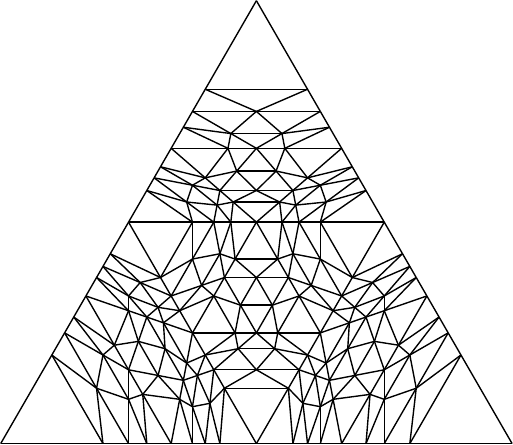}
\includegraphics[width=0.300000000000000\linewidth]{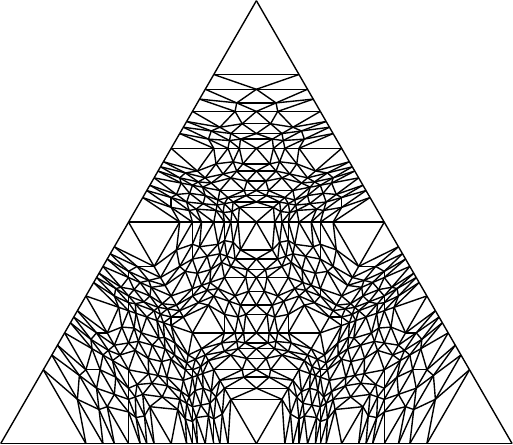}
\subsection{Density function}

The density function of the invariant measure of $f:\Delta\to\Delta$ for
the Reverse algorithm is \cite{arnoux_symmetric_2015}:
\[
\frac{1}{(1-x_1)(1-x_2)(1-x_3)}.
\]

\subsection{Invariant measure}
\includegraphics[width=0.800000000000000\linewidth]{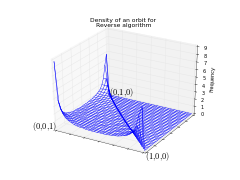}
\subsection{Natural extension}
Two sequences 
$(\mathbf{x}_{n+1})_{n\geq0}$ and 
$(\mathbf{a}_{n+1})_{n\geq0}$ defined such that
\[
    \mathbf{x}_{n+1}=M(\mathbf{x}_n)^{-1}\mathbf{x}_n
    \qquad\text{ and }\qquad
    \mathbf{a}_{n+1}=M(\mathbf{x}_n)^\top\mathbf{a}_n.
\]

\includegraphics[width=1\linewidth]{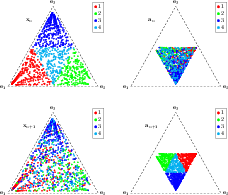}
\subsection{Lyapunov exponents}
(using 30 orbits of 
100000000 iterations each)\\
\begin{tabular}{lllll}
30 succesfull orbits & min & mean & max & std \\ \hline
$\theta_1$ & $0.4042$ & $0.4049$ & $0.4052$ & $0.00022$ \\
$\theta_2$ & $-0.10329$ & $-0.10319$ & $-0.10300$ & $0.000062$ \\
$1-\theta_2/\theta_1$ & $1.25472$ & $1.25487$ & $1.25498$ & $0.000055$ \\
\end{tabular}
\subsection{Substitutions}
\[
\begin{array}{lll}
\sigma_{1}=\left\{\begin{array}{l}
1 \mapsto 1\\
2 \mapsto 21\\
3 \mapsto 31
\end{array}\right.
&
\sigma_{2}=\left\{\begin{array}{l}
1 \mapsto 12\\
2 \mapsto 2\\
3 \mapsto 32
\end{array}\right.
&
\sigma_{3}=\left\{\begin{array}{l}
1 \mapsto 13\\
2 \mapsto 23\\
3 \mapsto 3
\end{array}\right.
\\
\sigma_{4}=\left\{\begin{array}{l}
1 \mapsto 23\\
2 \mapsto 31\\
3 \mapsto 12
\end{array}\right.
&
\end{array}
\]
\subsection{$S$-adic word example}
Using vector $v=\left(1, e, \pi\right)$:
\begin{align*}
w &=
\sigma_{4}
\sigma_{1}
\sigma_{1}
\sigma_{4}
\sigma_{3}
\sigma_{1}
\sigma_{1}
\sigma_{3}
\sigma_{3}
\sigma_{3}
\cdots(1)\\
& = 2331232331232312232323312323312323122323...
\end{align*}
Factor Complexity of $w$ is 
$(p_w(n))_{0\leq n \leq 20} =$
\[
(1, 3, 6, 9, 12, 14, 17, 20, 23, 26, 29, 32, 35, 38, 41, 44, 47, 50, 53, 56, 58)
\]
\subsection{Discrepancy}
ValueError: On input=[197, 2, 1], algorithm Reverse reaches non integer entries (0.5, 0.5, 1.5)
\subsection{Dual substitutions}
\[
\begin{array}{lll}
\sigma^*_{1}=\left\{\begin{array}{l}
1 \mapsto 123\\
2 \mapsto 2\\
3 \mapsto 3
\end{array}\right.
&
\sigma^*_{2}=\left\{\begin{array}{l}
1 \mapsto 1\\
2 \mapsto 231\\
3 \mapsto 3
\end{array}\right.
&
\sigma^*_{3}=\left\{\begin{array}{l}
1 \mapsto 1\\
2 \mapsto 2\\
3 \mapsto 312
\end{array}\right.
\\
\sigma^*_{4}=\left\{\begin{array}{l}
1 \mapsto 23\\
2 \mapsto 13\\
3 \mapsto 12
\end{array}\right.
&
\end{array}
\]
\subsection{E one star}
ValueError: The substitution (1->23, 2->1233, 3->1232) must be unimodular.
\newpage
\section{Cassaigne algorithm}
\subsection{Definition}
On $\Lambda=\mathbb{R}^3_+$, the map is
\cite{cassaigne_algorithme_2015}
\[
F (x_1,x_2,x_3) = 
\begin{cases}
    (x_1-x_3, x_3, x_2) & \mbox{if } x_1 > x_3\\
    (x_2, x_1, x_3-x_1) & \mbox{if } x_1 < x_3.
\end{cases}
\]
\subsection{Matrix Definition}
The partition of the cone is
$\Lambda=\cup_{\pi\in\mathcal{S}_3}\Lambda_\pi$ where
\begin{align*}
	\Lambda_1 &= \{(x_1,x_2,x_3)\in\Lambda\mid 
	x_1 > x_3\}, \\
    \Lambda_2 &= \{(x_1,x_2,x_3)\in\Lambda\mid 
	x_1 < x_3\}.
\end{align*}
The matrices are given by the rule
\[
    M(\mathbf{x}) = M_i
    \qquad\text{ if and only if }\qquad
    \mathbf{x}\in\Lambda_i.
\]
The map $F$ on $\Lambda$ and
the projective map $f$ on
$\Delta=\{\mathbf{x}\in\Lambda\mid\Vert\mathbf{x}\Vert_1=1\}$ are:
\[
    F(\mathbf{x}) = M(\mathbf{x})^{-1}\mathbf{x}
    \qquad\text{and}\qquad
    f(\mathbf{x}) = \frac{F(\mathbf{x})}{\Vert F(\mathbf{x})\Vert_1}.
\]

\subsection{Matrices}
\[
\begin{array}{lll}
M_{1}={\arraycolsep=2pt\left(\begin{array}{rrr}
1 & 1 & 0 \\
0 & 0 & 1 \\
0 & 1 & 0
\end{array}\right)}
&
M_{2}={\arraycolsep=2pt\left(\begin{array}{rrr}
0 & 1 & 0 \\
1 & 0 & 0 \\
0 & 1 & 1
\end{array}\right)}
&
\end{array}
\]
\subsection{Cylinders}
\includegraphics[width=0.300000000000000\linewidth]{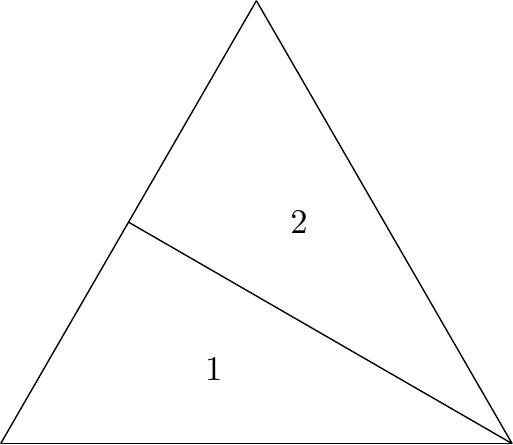}
\includegraphics[width=0.300000000000000\linewidth]{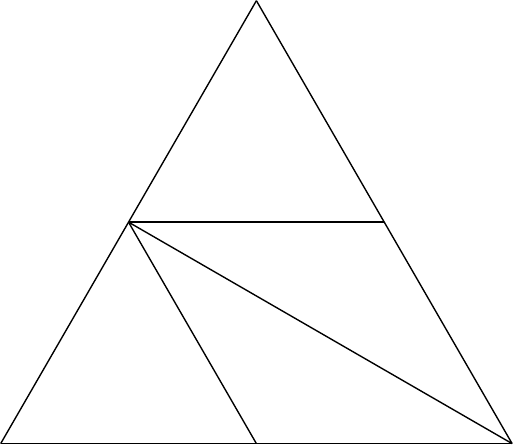}
\includegraphics[width=0.300000000000000\linewidth]{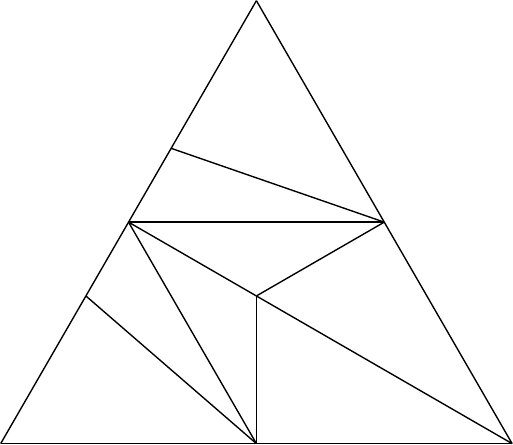}
\includegraphics[width=0.300000000000000\linewidth]{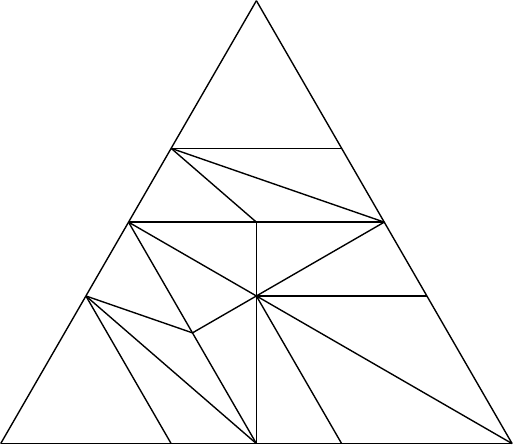}
\includegraphics[width=0.300000000000000\linewidth]{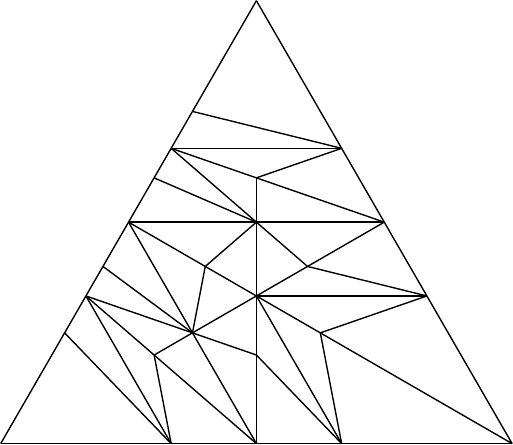}
\includegraphics[width=0.300000000000000\linewidth]{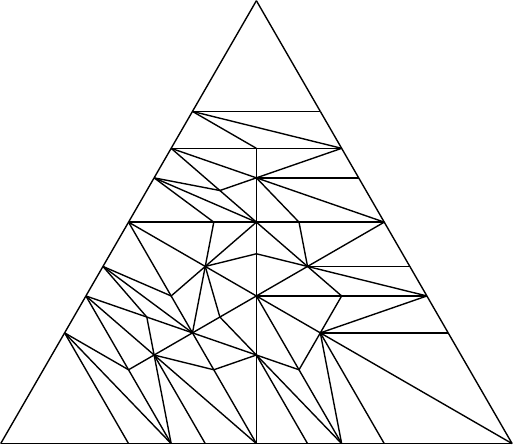}
\includegraphics[width=0.300000000000000\linewidth]{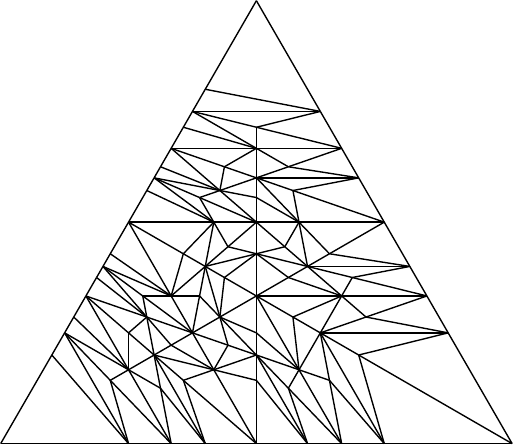}
\includegraphics[width=0.300000000000000\linewidth]{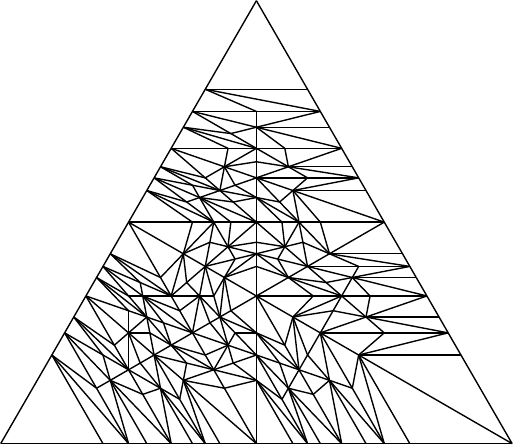}
\includegraphics[width=0.300000000000000\linewidth]{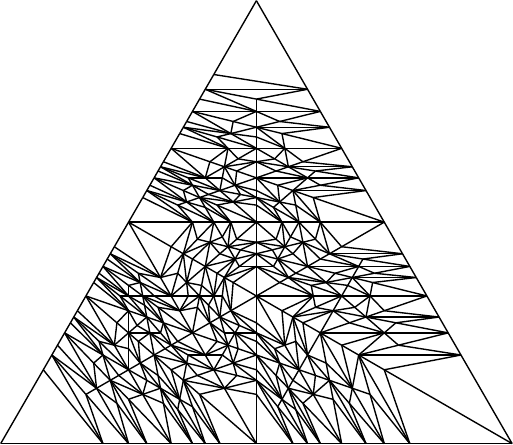}
\subsection{Density function}
The density function of the invariant measure of $f:\Delta\to\Delta$ for
the Cassaigne algorithm is
\cite{arnoux_symmetric_2015}
\[
\frac{1}{(1-x_1)(1-x_3)}.
\]

\subsection{Invariant measure}
\includegraphics[width=0.800000000000000\linewidth]{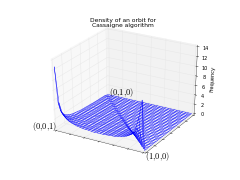}
\subsection{Natural extension}
Two sequences 
$(\mathbf{x}_{n+1})_{n\geq0}$ and 
$(\mathbf{a}_{n+1})_{n\geq0}$ defined such that
\[
    \mathbf{x}_{n+1}=M(\mathbf{x}_n)^{-1}\mathbf{x}_n
    \qquad\text{ and }\qquad
    \mathbf{a}_{n+1}=M(\mathbf{x}_n)^\top\mathbf{a}_n.
\]

\includegraphics[width=1\linewidth]{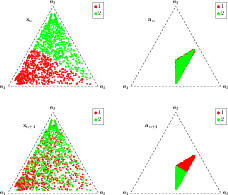}
\subsection{Lyapunov exponents}
(using 30 orbits of 
100000000 iterations each)\\
\begin{tabular}{lllll}
30 succesfull orbits & min & mean & max & std \\ \hline
$\theta_1$ & $0.1824$ & $0.1827$ & $0.1829$ & $0.00013$ \\
$\theta_2$ & $-0.07083$ & $-0.07072$ & $-0.07060$ & $0.000054$ \\
$1-\theta_2/\theta_1$ & $1.38698$ & $1.38712$ & $1.38725$ & $0.000070$ \\
\end{tabular}
\subsection{Substitutions}
\[
\begin{array}{lll}
\sigma_{1}=\left\{\begin{array}{l}
1 \mapsto 1\\
2 \mapsto 13\\
3 \mapsto 2
\end{array}\right.
&
\sigma_{2}=\left\{\begin{array}{l}
1 \mapsto 2\\
2 \mapsto 13\\
3 \mapsto 3
\end{array}\right.
&
\end{array}
\]
\subsection{$S$-adic word example}
Using vector $v=\left(1, e, \pi\right)$:
\begin{align*}
w &=
\sigma_{2}
\sigma_{1}
\sigma_{2}
\sigma_{1}
\sigma_{1}
\sigma_{1}
\sigma_{1}
\sigma_{2}
\sigma_{1}
\sigma_{1}
\cdots(1)\\
& = 2323213232323132323213232321323231323232...
\end{align*}
Factor Complexity of $w$ is 
$(p_w(n))_{0\leq n \leq 20} =$
\[
(1, 3, 5, 7, 9, 11, 13, 15, 17, 19, 21, 23, 25, 27, 29, 31, 33, 35, 37, 39, 41)
\]
\subsection{Discrepancy}
Discrepancy \cite{MR593979} for all 19701
 $S$-adic words with directions
$v\in\mathbb{N}^3_{>0}$
such that $v_1+v_2+v_3=200$:
\begin{center}
\includegraphics[width=0.600000000000000\linewidth]{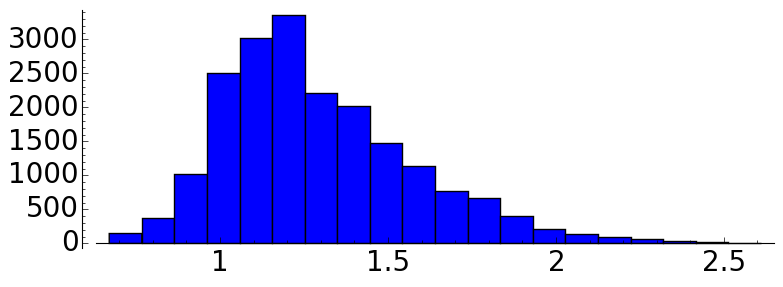}
\end{center}
\subsection{Dual substitutions}
\[
\begin{array}{lll}
\sigma^*_{1}=\left\{\begin{array}{l}
1 \mapsto 12\\
2 \mapsto 3\\
3 \mapsto 2
\end{array}\right.
&
\sigma^*_{2}=\left\{\begin{array}{l}
1 \mapsto 2\\
2 \mapsto 1\\
3 \mapsto 23
\end{array}\right.
&
\end{array}
\]
\subsection{E one star}
Using vector $v=\left(1, e, \pi\right)$, the 13-th 
iteration on the unit cube is:
\[
E_1^*(\sigma^*_{2})
E_1^*(\sigma^*_{1})
E_1^*(\sigma^*_{2})
E_1^*(\sigma^*_{1})
E_1^*(\sigma^*_{1})
\cdots
(\includegraphics[width=1em]{cube.pdf})=
\]
\begin{center}
\includegraphics[height=3cm]{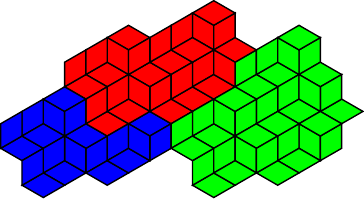}
\end{center}
\newpage

\end{multicols}

\footnotesize
\section{Comparison of Lyapunov exponents}
(30 orbits of 
1000000000 iterations each)
\begin{center}
\begin{tabular}{lllll}
Algorithm & \#Orbits & $\theta_1$ (std) & $\theta_2$ (std) & $1-\theta_2/\theta_1$ (std) \\ \hline
Arnoux-Rauzy-Poincar\'e & $30$ & 0.44290 (0.000083) & -0.17219 (0.000035) & 1.38879 (0.000017) \\
Selmer & $30$ & 0.18269 (0.000032) & -0.07072 (0.000013) & 1.38710 (0.000021) \\
Cassaigne & $30$ & 0.18268 (0.000041) & -0.07072 (0.000017) & 1.38709 (0.000028) \\
Brun & $30$ & 0.30449 (0.000049) & -0.11216 (0.000019) & 1.36833 (0.000015) \\
Reverse & $30$ & 0.40489 (0.000057) & -0.10320 (0.000015) & 1.25489 (0.000016) \\
Fully Subtractive & $26$ & 2.5e-9 (1.6e-9) & 9.3e-8 (1.5e-9) & -69. (87.) \\
Poincar\'e & $22$ & 6.9e-9 (4.8e-9) & 7.8e-8 (3.2e-8) & -24. (40.) \\
\end{tabular}
\end{center}

\begin{multicols}{3}
\setlength{\premulticols}{1pt}
\setlength{\postmulticols}{1pt}
\setlength{\multicolsep}{1pt}
\setlength{\columnsep}{2pt}
\newcommand{\note}[1]{\hfill\textrm{\textcolor{gray}{#1}}}
\newcommand{\args}[1]{\textit{\textcolor{blue}{#1}}}
\newcommand{\stdout}[1]{\textcolor{Sepia}{#1}}
\raggedright
\footnotesize

\section{Sage Code}
This section shows how to reproduce any of the results in these Cheat Sheets.
\subsection{Requirements}
The image and experimental results in these cheat sheets were created with the
following version of Sage \cite{sage}
\begin{verbatim}
$ sage -v
SageMath Version 6.10.beta3, Release Date: 2015-11-05
\end{verbatim}
and my optional Sage package \cite{labbe_slabbe_2015} which can be installed with:
\begin{verbatim}
$ sage -p http://www.slabbe.org/Sage/slabbe-0.2.spkg
\end{verbatim}
\subsection{Definition}
Define a Multidimensional Continued Fraction algorithm:
\begin{verbatim}
sage: from slabbe.mult_cont_frac import Brun
sage: algo = Brun()
\end{verbatim}
You may replace \texttt{Brun} above by any of the following:
\begin{verbatim}
Brun, Poincare, Selmer, FullySubtractive, 
ARP, Reverse, Cassaigne
\end{verbatim}
\subsection{Matrices}
\begin{verbatim}
sage: cocycle = algo.matrix_cocycle()
sage: cocycle.gens()
\end{verbatim}
\subsection{Cylinders}
\begin{verbatim}
sage: cocycle = algo.matrix_cocycle()
sage: t = cocycle.tikz_n_cylinders(3, scale=3)
sage: t.pdf()
\end{verbatim}
\subsection{Density function}
This section is hand written.
\subsection{Invariant measure}
\begin{verbatim}
sage: fig = algo.invariant_measure_wireframe_plot(
....:       n_iterations=10^6, ndivs=30, norm='1')
sage: fig.savefig('a.pdf')
\end{verbatim}
\subsection{Natural extension}
\begin{verbatim}
sage: t = algo.natural_extension_tikz(n_iterations=1200, 
....:          marksize=.8, group_size="2 by 2")
sage: t.png()
\end{verbatim}
\subsection{Lyapunov exponents}
The algorithm that computes Lyapunov exponents was provided to me
by Vincent Delecroix, in June 2013. I translated his C code into cython.
\begin{verbatim}
sage: from slabbe.lyapunov import lyapunov_table
sage: lyapunov_table(algo, n_orbits=30, n_iterations=10^7)
\end{verbatim}
\subsection{Substitutions}
\begin{verbatim}
sage: algo.substitutions()
\end{verbatim}
\subsection{$S$-adic word example}
\begin{verbatim}
sage: v = (1,e,pi)
sage: it = algo.coding_iterator(v)
sage: [next(it) for _ in range(10)]
sage: algo.s_adic_word(v)
sage: map(w[:10000].number_of_factors, range(21))  
\end{verbatim}
\subsection{Discrepancy}
\begin{verbatim}
sage: D = algo.discrepancy_statistics(length=20)
sage: histogram(D.values())
\end{verbatim}
\subsection{Dual substitutions}
\begin{verbatim}
sage: algo.dual_substitutions()
\end{verbatim}
\subsection{E one star}
\begin{verbatim}
sage: from slabbe import TikzPicture
sage: P = algo.e_one_star_patch(v=(1,e,pi), n=8)
sage: s = P.plot_tikz()
sage: TikzPicture(s).pdf()
\end{verbatim}
\subsection{Comparison of Lyapunov exponents}
\begin{verbatim}
sage: import slabbe.mult_cont_frac as mcf
sage: from slabbe.lyapunov import lyapunov_comparison_table
sage: algos = [mcf.Brun(), mcf.Selmer(), mcf.ARP(), 
....:          mcf.Reverse(), mcf.Cassaigne()]
sage: lyapunov_comparison_table(algos, n_orbits=30, 
....:                           n_iterations=10^7)
\end{verbatim}

\section*{Acknowledgments}

This work is part of the project ``Dynamique des algorithmes du pgcd : une
approche Algorithmique, Analytique, Arithmétique et Symbolique (Dyna3S)''
(ANR-13-BS02-0003)  supported by the Agence Nationale de la Recherche. The
author is supported by a postdoctoral Marie Curie fellowship (BeIPD-COFUND)
cofunded by the European Commission. I wish to thank Valérie Berthé, Pierre
Arnoux, Vincent Delecroix and Thierry Monteil for many discussions on the
experimental aspects of MCF algorithms.

\rule{0.3\linewidth}{0.25pt}
\scriptsize
\bibliographystyle{plain}
\bibliography{biblio}

\begin{thebibliography}{10}

\bibitem{arnoux_symmetric_2015}
Pierre Arnoux and Sébastien Labbé.
\newblock On some symmetric multidimensional continued fraction algorithms.
\newblock {\em arXiv:1508.07814}, August 2015.

\bibitem{2015_berthe_factor}
V.~Berthé and S.~Labbé.
\newblock Factor complexity of {S}-adic words generated by the
  arnoux–rauzy–poincaré algorithm.
\newblock {\em Advances in Applied Mathematics}, 63(0):90 -- 130, 2015.

\bibitem{MR0111735}
Viggo Brun.
\newblock Algorithmes euclidiens pour trois et quatre nombres.
\newblock In {\em Treizi\`eme congr\`es des math\`ematiciens scandinaves, tenu
  \`a {H}elsinki 18-23 ao\^ut 1957}, pages 45--64. Mercators Tryckeri,
  Helsinki, 1958.

\bibitem{cassaigne_algorithme_2015}
Julien Cassaigne.
\newblock Un algorithme de fractions continues de complexité linéaire.
\newblock October 2015.
\newblock DynA3S meeting, LIAFA, Paris, October 12th, 2015.

\bibitem{labbe_slabbe_2015}
Sébastien Labbé.
\newblock Sébastien labbé research code v0.2, \texttt{slabbe-0.2.spkg}.
\newblock \url{http://www.slabbe.org/Sage}, 2015.

\bibitem{MR1336331}
A.~Nogueira.
\newblock The three-dimensional {P}oincar\'e continued fraction algorithm.
\newblock {\em Israel J. Math.}, 90(1-3):373--401, 1995.

\bibitem{schweiger}
F.~Schweiger.
\newblock {\em Multidimensional Continued Fraction}.
\newblock Oxford Univ. Press, New York, 2000.

\bibitem{MR0130852}
Ernst~S. Selmer.
\newblock Continued fractions in several dimensions.
\newblock {\em Nordisk Nat. Tidskr.}, 9:37--43, 95, 1961.

\bibitem{sage}
William~A. Stein et~al.
\newblock {\em {S}age {M}athematics {S}oftware ({V}ersion 6.9)}.
\newblock The Sage Development Team, 2015.

\bibitem{MR593979}
R.~Tijdeman.
\newblock The chairman assignment problem.
\newblock {\em Discrete Math.}, 32(3):323--330, 1980.

\end{thebibliography}

\end{multicols}
\end{document}